\documentclass[11pt,fleqn]{article}
\usepackage{amssymb,amsthm,epsfig}

\def\cro{{\mbox {\sc cr}}}
\def\gcro{{\mbox {\sc gcr}}}
\def\rcro{{\mbox {\sc rcr}}}

\begin{document}
\title{Note on $k$-planar crossing numbers\thanks{Research on this paper was conducted at the workshop on \emph{Exact Crossing Numbers}, April 28--May 2, 2014, at the American Institute of Mathematics, Palo Alto, CA.}}

\author{
J\'anos Pach\thanks{Ecole Polytechnique F\'ed\'erale de Lausanne, Station 8, CH-1015 Lausanne, Switzerland and R\'enyi Institute of Mathematics, Hungarian Academy of Sciences, PO Box 127, H-1364, Budapest, Hungary. Email: \texttt{pach@cims.nyu.edu}. Partially supported the by
\textsc{SNF} grants 200020-144531 and 200021-137574.}
\and
L\'aszl\'o A. Sz\'ekely\thanks{Department of Mathematics, University of South Carolina, Columbia, SC, USA.
Email: \texttt{szekely@math.sc.edu}. Partially supported by the \textsc{NSF} grant DMS~1300547.}
   \and
Csaba D. T\'oth\thanks{Department of Mathematics, California State University Northridge, Los Angeles, CA, USA and Department of Computer Science, Tufts University, Medford, MA, USA. Email: \texttt{csaba.toth@csun.edu}. Partially supported by the \textsc{NSF} awards CCF~1422311 and CCF~1423615.}
\and
G\'eza T\'oth\thanks{R\'enyi Institute of Mathematics, Hungarian Academy of Sciences, PO Box 127, H-1364 Budapest, Hungary. Email: \texttt{toth.geza@renyi.mta.hu}. Partially supported by the
\textsc{OTKA} grant K-83767.}
}

\date{\it Dedicated to our colleague Ferran Hurtado (1951--2014)}
\maketitle


\begin{abstract} The {\em crossing number} $\cro(G)$ of a graph $G=(V,E)$ is
the smallest number of edge crossings over all drawings of $G$ in the
plane. For any $k\ge 1$, the {\em $k$-planar crossing number} of $G$,
$\cro_k(G)$, is defined as the minimum of
$\cro(G_0)+\cro(G_1)+\ldots+\cro(G_{k-1})$ over all graphs $G_0, G_1,\ldots,
G_{k-1}$ with $\cup_{i=0}^{k-1}G_i=G$.
It is shown that for every $k\ge 1$, we have $\cro_k(G)\le
\left(\frac{2}{k^2}-\frac1{k^3}\right)\cro(G)$. This bound does not remain
true if we replace the constant $\frac{2}{k^2}-\frac1{k^3}$ by any number
smaller than $\frac1{k^2}$. Some of the results extend to the rectilinear
variants of the $k$-planar crossing number.
\end{abstract}

\section{Introduction}

Selfridge (see~\cite{Ha61}) noticed that by Euler's polyhedral formula
$K_{11}$, the complete graph on $11$ vertices, cannot be written as the union
of two planar graphs. Later Battle, Harary, and Kodama~\cite{BaHK62} and
independently Tutte~\cite{Tu63a} proved that the same is true for $K_9$, but
not for $K_8$. This led Tutte~\cite{Tu63b} to introduce a new parameter, the
\emph{thickness} of a graph $G$, which is the minimum number of planar graphs
that $G$ can be decomposed into. The notion turned out to be relevant for VLSI
chip design, where it corresponds to the number of layers required for
realizing a network so that there is no crossing within a layer. Consult
Mutzel, Odenthal, and Scharbrodt~\cite{MuOS98} for a survey. If the thickness
of $G$ is at most $2$, $G$ is called \emph{biplanar}. Mansfield proved that it
is an NP-complete problem to decide whether a graph is biplanar; see~\cite{beineke,Ma83}.

\smallskip

A \emph{drawing} of a graph $G=(V,E)$ is a planar representation of $G$ such
that every vertex $v\in V$ corresponds to a point of the plane and every edge
$uv\in E$ is represented by a simple continuous curve between the points
corresponding to $u$ and $v$, which does not pass through any point
representing a vertex of $G$. We always assume for simplicity
that (1) no two curves share infinitely many points, (2) no two curves are
tangent to each other, and (3) no three curves pass through the same
point. The \emph{crossing number} of $G$ is defined as the minimum number of
edge crossings in a drawing of $G$, and is denoted by $\cro(G)$. For surveys,
see~\cite{schaefer,success}. Clearly, $G$ is planar if and only if $\cro(G)=0$.

The \emph{biplanar crossing number}, $\cro_2(G)$, of $G$ was defined by
Owens~\cite{owens} as the minimum sum of the crossing numbers of two graphs,
$G_0$ and $G_1$, whose union is $G$. For the VLSI applications, we imagine
that $G_0$ and $G_1$ are drawn (realized) in different planes.  If $G$ is
biplanar, its biplanar crossing number is $0$. The biplanar crossing number of
random graphs was studied by Spencer~\cite{spencer}. Czabarka, S\'ykora,
Sz\'ekely, and Vr\v to  \cite{bipII} proved that for every
graph $G$, we have
$$\cro_2(G)\leq \frac{3}{8}\cro(G).$$
They also showed~\cite{bipsurvey} that this inequality does not remain true
if the constant $\frac38=0.375$ is replaced by anything less than $\frac{8}{119}\approx 0.067$.

\smallskip

Shahrokhi et al.~\cite{kplanar} extended the notion of biplanar crossing
number as follows. For any positive integer $k\ge 1$, define the
\emph{$k$-planar crossing number} of $G$ as the minimum of
$\cro(G_0)+\cro(G_1)+\ldots+\cro(G_{k-1})$, where the minimum is taken over
all graphs $G_0, G_1,\ldots, G_{k-1}$ whose union is $G$, that is,
$\cup_{i=0}^{k-1}E(G_i)=E(G)$. This number is denoted by
$\cro_k(G)$. Obviously, $\cro_1(G)=\cro(G)$ and we have
$\cro_i(G)\ge\cro_{i+1}(G)$ for all $i\in\mathbb{N}$ and every graph $G$.

\smallskip

In the present note, we investigate the relationship between the $k$-planar
crossing number and the (ordinary) crossing number of a graph.
For every $k\ge 1$, let
$$\alpha_k= \sup {\cro_k(G)\over \cro(G)},$$
where the supremum is taken over all \emph{nonplanar} graphs $G$.
The above mentioned results yield $0.067<\alpha_2\le \frac38=0.375$.
The next theorem implies that $\alpha_k= \Theta(k^{-2})$.

\medskip
\noindent {\bf Theorem.} {\em For every positive integer $k$, we have}
$${1\over k^2}\le\alpha_k\le {2\over k^2}-{1\over k^3}.$$
\smallskip

\section{Proof of Theorem} \label{proofof}

\paragraph{Upper bound.} First we prove the upper bound.
Let $G$ be a graph with vertex set $V(G)$, edge set $E(G)$, and fix an optimal
drawing of $G$ in the plane with precisely
$\cro(G)$ crossings. We describe a randomized procedure to partition (the edge set of) $G$ into
$k$ subgraphs $G_0,\ldots, G_{k-1}$ such that the expected value of the sum of their crossing numbers
is at most $(\frac2{k^2}-\frac1{k^3})\cro(G)$. We think of each $G_i$ as a graph drawn independently
so that edges of different subgraphs do not cross.

\smallskip

The idea of the proof is the following. We start by randomly partitioning the
vertex set of $G$ into $k$ roughly equal classes. We associate with each class
a vertex of a complete graph $K_k$. We consider a factorization of $K_k$ into
maximal matchings and then use these matchings to divide $E(G)$ into $k$ classes,
$G_0,\ldots, G_{k-1}$. It will follow from the definition that every $G_i$ can
be drawn independently in such a way that no two edges that correspond
to distinct edges of the underlying matching of $K_k$ will cross.

\smallskip

Let the vertex set of $G$ be $V=V(G)=\{ 1, 2, \ldots , n\}$. Assign
independent random variables $\xi_v$ to the vertices $v\in V$ such that each
$\xi_v$ takes each of the values $0, 1, \ldots, k-1$ with probability $1/k$.

For every $i$($0\le i<k$), let $V_i=\{ v\in V\ |\ \xi_v=i\}$, and define a
subgraph $G_i$ as follows. Let $V(G_i)=V$ and let the edge set $E(G_i)$ of
$G_i$ consist of all edges $uv\in E(G)$ for which
$$\xi_u+\xi_v\equiv i\bmod k.$$
Obviously, we have $\cup_{i=0}^{k-1}E(G_i)=E(G)$.

\smallskip

We define the {\em type} of an edge $uv$ to be the unordered pair $(\xi_u, \xi_v)$.
For each $i\; (0\le i<k)$, first we draw $G_i$ in the $i$th plane as it was
drawn in the original drawing of $G$. Notice that for every index $g$, there
is precisely one index $h=h(g)$ such that $G_i$ has an edge connecting a vertex in
$V_g$ to a vertex in $V_h$. Thus, every connected component of $G_i$ consists
of edges of the same type. In the $i$th plane, we can translate the connected
components of $G_i$ sufficiently far from each other so that no two edges of
different types intersect, and during the procedure no new crossings are
introduced.

\smallskip

Calculate the expected value of the total number of crossings in the resulting
drawing of $G_i$ over all $i\; (0\le i<k)$. Every crossing arises from a
crossing between two edges in the original drawing of $G$.
Consider two edges $uv, u'v'\in E(G)$ that cross each other in the original drawing.
A crossing between these edges will be present in the final drawing of one of
the $G_i$s if and only $uv$ and $u'v'$ are of the same type. For every index $g$,
this happens with probability
${\mbox{Pr}}[{\mbox{type}}(uv)=(g, g)]={1\over k^2}$. For distinct indices
$g$ and $h\; (g\neq h)$, we have
${\mbox{Pr}}[{\mbox{type}}(uv)=(g, h)]={2\over k^2}$.

Summing over all possible pairs of types, we obtain
$${\mbox{Pr}}[{\mbox{type}}(uv)={\mbox{type}}(u'v')]
={k\choose 2}\cdot{2\over k^2}\cdot {2\over k^2}+
k\cdot{1\over k^2}\cdot {1\over k^2}
={2\over k^2}-{1\over k^3}.$$
Consequently, the expected value of the total number of crossings in the resulting
drawings of all $G_i$s is $({2\over k^2}-{1\over k^3})\cro(G)$.
Hence, there exists a partition of (the edges of) $G$ into $G_0,\ldots ,G_{k-1}$
where
$$\cro(G_0)+\ldots+\cro(G_{k-1})\le \left({2\over k^2}-{1\over k^3}\right)\cro(G).$$
This completes the proof of the upper bound in the Theorem.
\medskip

\paragraph{Lower bound.} Next we establish the lower bound.
For two functions $f(n)$ and $g(n)$, we write $f(n)\ll g(n)$, if $\lim_{n\rightarrow\infty}\frac{f(n)}{g(n)}=0$.
Let $\kappa(n,e)$ denote the minimum crossing number of a
graph $G$ with $n$ vertices and at least $e$ edges. That is,
$$\kappa(n,e)=\min_{\begin{array}{cc}
|V(G)|=n\\ |E(G)|\ge e \end{array}}\cro(G).$$

\medskip

It was shown in \cite{PST} that there
exists a positive constant $K$ such that if $n\ll e\ll n^2$, the limit
$$\lim_{n\rightarrow\infty}\kappa(n,e){n^2\over e^3}$$
exists and is equal to $K$.
The constant $K>0$ is called the \emph{midrange crossing constant}. 
The best known bounds for $K$ are $0.034\le K\le 0.09$; see \cite{A15,PRTT,PT}.
This result can be rephrased as follows.

\medskip

\noindent {\bf Lemma.} {\em For every $\varepsilon\; (0<\varepsilon<1)$,
there exists a constant $N=N_{\varepsilon}$ satisfying the following condition. For every positive integers $n$ and $e$ with $\min(n, \frac{e}{n}, \frac{n^2}{e})\ge N$, we have $\kappa(n,e)>(K-\varepsilon)\frac{e^3}{n^2}$, and there is a graph $G$ with $n$ vertices and $e$ edges such that $\cro(G)< (K+\varepsilon )\frac{e^3}{n^2}$.}

\medskip

Let $\varepsilon>0$ be fixed, let
$$\min\left(n, \frac{e}{n}, \frac{n^2}{e}\right) > \frac{k}{\varepsilon}N_{\varepsilon},$$ and let
$G$ be a graph with  $n$ vertices and $e$ edges such that
$\cro(G)< (K+\varepsilon )\frac{e^3}{n^2}$.
Decompose $G$ into $k$ graphs $G=G_0\cup G_1,\cdots \cup G_{k-1}$
such that $\cro(G_0)+\cro(G_1)+\cdots +\cro(G_{k-1})=\cro_k(G)$.
For simplicity, write $e_i$ for $|E(G_i)|$.
\smallskip

We may assume, without loss of generality, that there is an integer $t\; (0<t\le k)$ such that $e_i\ge \frac{\varepsilon}{k}e$ for $i=0, 1,\ldots, t-1$, and $e_i< \frac{\varepsilon}{k}e$ for $i=t, t+1,\ldots, k-1$.

For every $i<t$, we have $\min(n, \frac{e_i}{n}, \frac{n^2}{e_i}) > N_{\varepsilon}$,
so we can apply the Lemma to conclude that $\cro(G_i)\ge (K-\varepsilon)\frac{e_i^3}{n^2}$.
Using that $\sum_{i=t}^{k-1}e_i\le \varepsilon e$, we have
$\sum_{i=0}^{t-1}e_i\ge (1-\varepsilon)e$.

Hence, Jensen's inequality yields
\begin{eqnarray*}
%
\cro_k(G)& \ge &\sum_{i=0}^{t-1}\cro(G_i)\geq \sum_{i=0}^{t-1}(K-\varepsilon)\frac{e_i^3}{n^2}\\
&\geq &t(K-\varepsilon )\cdot \frac{((1-\varepsilon )e/t)^3}{n^2}>
\frac{(1-3\varepsilon )(K-\varepsilon )}{k^2}\cdot \frac{e^3}{n^2}.
\end{eqnarray*}

Using that $\cro(G)< (K+\varepsilon )\frac{e^3}{n^2}$, the last inequality implies
$$\frac{\cro_k(G)}{\cro(G)}\geq (1-3\varepsilon)\frac{K-\varepsilon}{K+\varepsilon}\cdot \frac{1}{k^2}.$$

As $\varepsilon\rightarrow 0$, the lower bound in the Theorem follows.

\section{Rectilinear Variants}\label{conclusion}

\paragraph{Rectilinear $k$-planar crossing numbers.}
The \emph{rectilinear crossing number}, $\rcro(G)$, of a graph $G$ is the minimum number of crossings
over all \emph{straight-line} drawings of $G$, in which the edges are represented by line segments. Obviously,
we have $\cro(G)\leq \rcro(G)$ for every graph $G$. For every $t\geq 4$, Bienstock and Dean~\cite{BD93} constructed families of graphs whose crossing number is at most $t$ and whose rectilinear crossing number is unbounded.

Similarly to $\cro_k(G)$, we define the \emph{rectilinear $k$-planar crossing number} of a graph $G$, denoted $\rcro_k(G)$, as the minimum of $\rcro(G_0)+\rcro(G_1)+\ldots+\rcro(G_{k-1})$, where the minimum is taken over
all graphs $G_0, G_1,\ldots, G_{k-1}$ whose union is $G$. It is clear that $\cro_k(G)\leq \rcro_k(G)$ for
every $k\in \mathbb{N}$. However, we do not know of any graph $G$ where $\cro_k(G)<\rcro_k(G)$ and $k\geq 2$.

The analogue of $\alpha_k$ for every $k\in \mathbb{N}$ is
$$\beta_k= \sup {\rcro_k(G)\over \rcro(G)},$$
where the supremum is taken over all \emph{nonplanar} graphs $G$. The proof of our main theorem carries
over verbatim to this variant, and yields
$${1\over k^2}\le\beta_k\le {2\over k^2}-{1\over k^3}.$$
Specifically, the upper bound starts from a fixed straight-line drawing of $G$ with exactly $\rcro(G)$ crossings. Our randomized procedure decomposes $G$ into $k$ graphs $G_0,\dots , G_{k-1}$, each of which consists of $k$ vertex-disjoint subgraphs induced by the $k$ edge types. These $k^2$ subgraphs can be translated independently to avoid any crossings between edges of different subgraphs, but maintain a straight-line drawing for each. The lower bound relies on the existence of a midrange crossing constant $\overline{K}>0$ for the \emph{rectilinear} crossing number, which is established by the argument in~\cite{PST} even though the constants $K$ and $\overline{K}$ are not necessarily the same.

\paragraph{Geometric $k$-planar crossing numbers.}
The \emph{geometric thickness} of a graph $G$, introduced by Kainen~\cite{Kai73}, is the smallest positive integer $k$ such that $G$ admits a $k$-edge-coloring \emph{and} a straight-line drawing in which edges of the same color do not cross. The color classes define a decomposition of $G$ into $k$ planar graphs $G_0,\ldots , G_{k-1}$ each of which admits a crossing-free straight-line drawing in such a way that corresponding vertices are represented by the same point in the plane. A straight-line drawing of a graph $G$ is called \emph{biplane} if $G$ admits a 2-edge-coloring such that no two edges of the same color cross in this drawing; see~\cite{GHK15}. Eppstein~\cite{Epp04} constructed graphs with thickness 3 and geometric thickness at least $t$ for  every $t>0$. Determining the geometric thickness of a graph is also an NP-hard problem~\cite{DGM13}.

The geometric thickness motivates the following variant of the $k$-planar crossing number.
The \emph{geometric $k$-planar crossing number} of a graph $G$, denoted $\gcro_k(G)$, is the minimum number of crossings between edges of the same color over all $k$-edge-colorings of $G$ and all straight-line drawings of $G$. It is clear that $\cro_k(G)\leq \rcro_k(G)\leq \gcro_k(G)$ for every graph $G$ and every $k\in \mathbb{N}$.

The analogue of $\alpha_k$ for every $k\in \mathbb{N}$ is
$$\gamma_k= \sup {\gcro_k(G)\over \rcro(G)},$$
where the supremum is taken over all \emph{nonplanar} graphs $G$.
The lower bound of our main theorem carries over verbatim to this variant, since it relies on density results, namely the (rectilinear) midrange crossing number. However, the upper bound argument does not extend to this variant. Our randomized procedure partitions the edge set $E(G)$ into $k$ color classes $E(G_0),\ldots ,E(G_{k-1})$, and crossings between edges of different colors do not count. But each color class consists of edges of up to $k$ different types, and the crossings between edges of the same color and different types cannot be eliminated. A weaker upper bound easily follows from a uniform random $k$-coloring of the edges, and yields
$${1\over k^2}\le\gamma_k\le {1\over k}.$$

\end{document}